\documentclass[12pt]{article}
\usepackage{amsmath}
\usepackage{amssymb}

\usepackage{graphicx}

\newcommand{\eps}{\varepsilon}
\newcommand{\Z}{{\mathbb Z}}

\newcommand{\G}{{\mathcal G}}
\newcommand{\Sph}{{\mathbb S}}

\newcommand{\MM}{{\mathcal M}}
\newcommand{\R}{{\mathbb R}}

\newcommand{\FFF}{{\mathcal F}}
\newcommand{\KK}{{\mathcal K}}

\newcommand{\SSS}{\mathcal{S}}
\renewcommand{\phi}{\varphi}

\newcommand{\om}{\text{\boldmath ${\omega}$}}
\newcommand{\s}{\text{\boldmath ${\sigma}$}}
\newcommand{\tb}{\text{\boldmath ${\theta}$}}

\newcommand{\EA}{{\mathbb E}}
\newcommand{\PA}{{\mathbb P}}
\newcommand{\EE}{{\mathbf E}}
\newcommand{\PP}{{\mathbf P}}

\newcommand{\Ps}{{\mathtt P}_{\!\text{\boldmath ${\sigma}$}}}
\newcommand{\Es}{{\mathtt E}_{\text{\boldmath ${\sigma}$}}}

\newcommand{\Po}{P_{\!\text{\boldmath ${\omega}$}}}
\newcommand{\Pt}{P_{\!\text{\boldmath ${\text{\boldmath ${\theta}$}}$}}}

\newcommand{\1}[1]{{\mathbf 1}{\{#1\}}}

\allowdisplaybreaks

\newtheorem{theo}{Theorem}[section]
\newtheorem{lmm}{Lemma}[section]
\newtheorem{df}{Definition}[section]
\newtheorem{prop}{Proposition}[section]

\newtheorem{rmk}{Remark}[section]

\newcommand{\supp}{\mathop{\rm supp}}
\newcommand{\qed}{\hfill$\Box$}

\title{Survival time of random walk in random environment among soft obstacles}
\author{Nina Gantert$^{1}$ \and
Serguei~Popov$^{2}$ \and 
Marina Vachkovskaia$^{3}$}

\begin{document}

\maketitle
{\footnotesize 
\smallskip
\noindent $^{~1}$CeNos Center for Nonlinear Science and Institut f\"{u}r Mathematische Statistik,
Fachbe\-reich Mathematik und Informatik,
Ein\-stein\-strasse 62, 48149 M\"{u}nster,
Germany\\ 
\noindent e-mail: \texttt{gantert@math.uni-muenster.de}

\smallskip
\noindent $^{~2}$Instituto de Matem{\'a}tica e Estat{\'\i}stica,
Universidade de S{\~a}o Paulo, rua do Mat{\~a}o 1010, CEP 05508--090,
S{\~a}o Paulo SP, Brasil\\
\noindent e-mail: \texttt{popov@ime.usp.br}, 
\noindent url: \texttt{www.ime.usp.br/$\sim$popov}

\smallskip
\noindent $^{~3}$Departamento de Estat\'\i stica, Instituto de Matem\'atica,
Estat\'\i stica e Computa\c{c}\~ao Cien\-t\'\i\-{}fica,
Universidade de Campinas,
Caixa Postal 6065, CEP 13083--970, Campinas SP, Brasil\\
\noindent e-mail: \texttt{marinav@ime.unicamp.br},
\noindent url: \texttt{www.ime.unicamp.br/$\sim$marinav}
}

\begin{abstract}
We consider a Random Walk in Random Environment (RWRE) moving in an i.i.d.\ random field of obstacles. When the particle hits an obstacle, it disappears with 
a positive probability. We obtain quenched and annealed bounds on the tails of the survival time in the general $d$-dimensional case. We then consider
a simplified one-dimensional model (where transition probabilities and obstacles are independent and the RWRE only moves to neighbour sites), and obtain finer results for the tail of the survival time. In addition, we study also the ``mixed" probability 
measures (quenched with respect to the obstacles and
annealed with respect to the transition probabilities and vice-versa) and give results for tails of the survival time with respect to these probability measures. Further, we apply the same methods to obtain bounds for the tails of hitting times of Branching Random 
Walks in Random Environment (BRWRE).
\\[.3cm]{\bf Keywords:} confinement of RWRE, survival time, quenched and annealed tails, nestling RWRE, branching 
random walks in random environment.
\\[.3cm]{\bf AMS 2000 subject classifications:} 60K37
\\[.3cm] Submitted: June 10, 2008, accepted: January 20, 2009
\end{abstract}

\section{Introduction and main results}
\label{s_intro}
Random walk and Brownian motion among random obstacles have been investigated intensively in the last three decades. For an introduction to the subject, its connections with other areas and an exposition of the techniques used, we refer to the book \cite{SZ}. Usually, one distinguishes \emph{hard obstacles}, where the particle is killed upon hitting them, and \emph{soft obstacles} where the particle is only killed with a certain probability. A typical model treated extensively in~\cite{SZ} is Brownian motion in a Poissonian field of obstacles. The following questions arise for this model: what is the asymptotic behaviour of the survival time? What is the best strategy of survival, i.e.\ what is the conditioned behaviour of the particle, given that it has survived until time~$n$? An important role in answering these questions has been played by the concept of ``pockets of low local eigenvalues'' (again, we refer to~\cite{SZ} for explanations). A key distinction in random media is the difference between the \emph{quenched} probability measure (where one fixes the environment) and the \emph{annealed} probability measure (where one averages over the environment).

In this paper, we are considering a discrete model with soft obstacles where there are two sources of randomness in the environment: the 
particle has random transition probabilities (which are assigned to the sites of the lattice in an i.i.d.\ way and then fixed for all times) and the obstacles are also placed randomly on the lattice and then their positions remain unchanged. We investigate the tails of the survival time. Similar questions have been asked in \cite{A1} for simple random walk. The ``pockets of low local eigenvalues'' are in our case ``traps free of obstacles'': these are regions without obstacles, where the transition probabilities are such that the particle tends to spend a long time there before escaping. These regions are responsible for the possibility of large survival time. We assume that the environments (transition probabilities and obstacles)  in all sites are independent and obtain quenched and annealed bounds on the tails of the survival time in the general $d$-dimensional case. We then consider a simplified one-dimensional model (where transition probabilities and obstacles are independent and the RWRE only moves to neighbour sites), and obtain finer results for the tail of the survival time. Having two sources of randomness in the environment, we study also the ``mixed" probability 
measures (quenched with respect to the obstacles and
annealed with respect to the transition probabilities and vice-versa) and give results for the tails of the survival time with respect to these probability measures. Further, we develop the analogy with the branching random walks in random environment (BRWRE) \cite{CP07-1,CP07-2}, and provide quenched and annealed bounds for hitting times in the BRWRE model.

Now we define the model formally.
Denote by $e_1, \ldots, e_d$ the coordinate vectors, and let $\|\cdot\|_1$ and $\|\cdot\|_2$
stand for the~$L^1$ and~$L^2$ norms in~$\Z^d$ respectively.
The \emph{environment} consists of the set
of transition probabilities $\om=(\omega_x(y),\; x,y\in \Z^d)$,
and the set of variables indicating the locations of the obstacles
$\tb=(\theta_x,\; x\in \Z^d)$, where $\theta_x=\1{\text{there is an obstacle in $x$}}$.
Let us denote by $\sigma_x=(\omega_x(\cdot), \theta_x)$ the environment in~$x\in\Z^d$,
and $\s=(\om, \tb)=(\sigma_x, \;x\in \Z^d)$ stands for the (global) environment.
We suppose that jumps are uniformly bounded by some constant~${\widetilde m}$, 
which means that $\omega_x(y)=0$ if $\|x-y\|_1> {\widetilde m}$.
Let us denote by~$\MM$ the environment space where the $\sigma_x$ are defined, i.e.\ 
\begin{align*}
\MM =& \Big\{\big((a(y))_{y\in\Z^d}, b\big): a(y)\ge 0, \forall y\in\Z^d,  
        \\      
  &\qquad\qquad\quad\;\sum_{y\in \Z^d} a(y)=1 =\!\!\! 
 \sum_{y:\|y\|_1\leq {\widetilde m} } a(y), \quad b\in\{0,1\}\Big\}.
\end{align*}
We assume that $(\sigma_x, \;x\in \Z^d)$ is a
collection of i.i.d.\ random variables. We denote by~$\PP$ the corresponding product measure, and by~$\EE$ its expectation.
In some cases these assumptions may be relaxed, see Remark~\ref{moge} below.

Let 
\[
p = \PP[\theta_0 = 1].
\]
We always assume $0 < p < 1$. 

\smallskip

Having fixed the realization of the random environment~$\s$, we now 
define the random walk~$\xi$ and the random time~$\tau$ as follows.
The discrete time random walk $\xi_n$ starts from some $z_0\in\Z^d$ and moves according to the 
 transition probabilities
\[
\Ps^{z_0}[\xi_{n+1}=x+y\mid \xi_n=x]=\omega_x(y).
\]
Here $\Ps^{z_0}$ stands for the so-called quenched probability (i.e., with fixed environment~$\s$)
and we denote by $\Es^{z_0}$ the corresponding expectation. Usually, we shall 
assume that the random walk starts from the origin, so that $z_0=0$; in this
case we use the simplified notations $\Ps, \Es$.

Fix $r \in (0,1)$ and let $Z_1, Z_2, \ldots$ be a sequence of i.i.d.\
 Bernoulli random variables with $\Ps[Z_i=1]=r$.
Denote by $\Theta=\{x\in \Z^d:\; \theta_x=1\}$ the set of sites where the obstacles are placed.
Let 
\[
\tau=\min\{n: \; \xi_n\in \Theta, Z_n=1\}.
\]
Intuitively, when the RWRE hits an obstacle, it ``disappears'' with probability~$r$, 
and~$\tau$ is the survival time of the particle.

We shall also consider the annealed probability law $\PA^{z_0}=\PP\times \Ps^{z_0}[\cdot]$, and the 
corresponding expectation $\EA^{z_0}=\EE \Es^{z_0}$.
Again, when the random walk starts from the origin, we use the simplified
notations $\PA,\EA$.

Throughout this paper we suppose that that
the environment~$\s$ satisfies the following two conditions.

\medskip
\noindent
{\bf Condition E.} There exists $\eps_0>0$ such that $\omega_x (e)\ge \eps_0$
for all $e\in \{\pm e_i, \; i=1, \ldots, d\}$, $\PP$-a.s.

\medskip

For $(\omega, \theta) \in \MM$, let
\begin{equation}\label{Deltadef}
\Delta_\omega = \sum\limits_y y\omega(y)
\end{equation}
be the drift of $\omega$.

\medskip
\noindent
{\bf Condition N.} 
We have 
\[
\PP[\theta_0=0, \Delta_\omega \cdot a > 0] > 0\quad \text{for all }a \in \R^d\setminus\{0\}.
\]

\medskip

Condition~E is a natural (uniform) ellipticity condition; Condition~N is a standard condition for RWRE and assures that the environment has
``traps'' (i.e., pieces of the environment free from obstacles from where it takes a long time to 
escape). 
Let us emphasize that in this paper the term ``trap'' does not refer to the disappearing of the particle; on the contrary, 
by ``the particle is trapped in some place''
we usually mean that the particle stays alive but is likely to remain in that
place for a long time. 
Observe that Condition~N
 implies that the RWRE is (strictly) \emph{nestling}, i.e.\ 
the origin is in the interior of the convex hull of the support of~$\Delta_\omega$.

Our goal is to study the quenched and annealed tails of the distribution
of~$\tau$: $\Ps[\tau>n]$ and $\PA[\tau>n]$.

First, we formulate the results on the tails of~$\tau$ in the $d$-dimensional case under the above assumptions:
\begin{theo}
\label{t_ann}
For all~$d\geq 1$ there exist $K_{i}^a(d)>0$, $i=1,2$, such that for all~$n$
\begin{equation}
\label{up_ann}
 \PA[\tau>n] \leq \exp(-K_1^a(d) \ln^d n),
\end{equation}
and
\begin{equation}
\label{low_ann}
 \PA[\tau>n] \geq \exp(-K_2^a(d) \ln^d n).
\end{equation}
\end{theo}
\begin{theo}
\label{t_qu}
For $d\ge 1$ there exist $K_i^q(d)>0$, $i=1,2,3,4$
(also with the property $K_j^q(1)<1$ for $j=2,4$),  
such that for $\PP$-almost all~$\s$ there exists
$n_0(\s)$ such that for all~$n\geq n_0(\s)$ we have
\begin{equation}
\label{up_qu_d}
 \Ps[\tau>n] \leq \exp\big(-K_1^q(d) n \exp(-K^q_2(d) \ln^{1/d} n)\big),
\end{equation}
and 
\begin{equation}
\label{low_qu_d}
 \Ps[\tau>n] \geq \exp\big(-K_3^q(d) n \exp(-K^q_4(d) \ln^{1/d} n)\big).
\end{equation}
\end{theo}

In fact (as it will be discussed in Section~\ref{s_proofs_d}) the original motivation for the model 
of this paper came from the study of the hitting times for branching random walks in random environment (BRWRE),
see~\cite{CP07-1}. The above Theorem~\ref{t_ann} has direct counterparts in~\cite{CP07-1}, namely,
Theorems~1.8 and~1.9. However, the problem of finding upper and lower bounds on the quenched tails of the hitting
times was left open in that paper (except for the case $d=1$). 
Now, Theorem~\ref{t_qu} of the present paper allows us to obtain the ana\-lo\-gous bounds also for the model of~\cite{CP07-1}. 
The model of~\cite{CP07-1} can be described as follows. 
Particles live in~$\Z^d$ and evolve in discrete time.
At each time, every particle in a site is substituted by (possibly more than one but at least one)
offspring which are placed in neighbour sites, independently of the other particles.
The rules of offspring generation 
 (similarly to the notation of the present paper, they are given by $\omega_x$ at site $x$) 
depend only on the location of the particle.
Similarly to the situation of this paper,
the collection $\om$ of those rules (the \emph{environment}) is itself random,
it is chosen in an i.i.d.\ way before starting the process, and then it is kept
fixed during all the subsequent evolution of the particle system. 
We denote by $\omega$ a generic element of the set of all possible environments 
at a given point,
and we distinguish $\omega$ with branching (the particle can be replaced with several particles) and $\omega$ without branching (the particle can only be replaced with exactly one particle).
The BRWRE is called recurrent if for almost all environments in the process (starting with one particle at the origin), all sites
are visited (by some particle) infinitely often a.s.
Using the notations of~\cite{CP07-1} (in particular, ${\mathtt P}_{\text{\boldmath ${\omega}$}}$
stands for the quenched probability law of the BRWRE),
we have the following
\begin{prop}
\label{brwre}
 Suppose that the BRWRE
is recurrent and uniformly elliptic. 
Let $T(0,x_0)$ be the hitting time of~$x_0$ for the BRWRE starting from~$0$. Then, there exist ${\widehat K}_i^q(d)>0$, $i=1,2,3,4$, such that for almost all environments $\om$,
there exists
$n_0(\om)$ such that for all~$n\geq n_0(\om)$ we have
\begin{equation}
\label{brwre_up}
 {\mathtt P}_{\text{\!\boldmath ${\omega}$}}[T(0,x_0)>n] \leq \exp\big(-{\widehat K}_1^q(d) n 
              \exp(-{\widehat K}^q_2(d) \ln^{1/d} n)\big).
\end{equation}
Now, let~$\G$ be the set of all $\omega$ without branching. Suppose that it has positive probability and
the origin belongs to the 
interior of the convex hull of $\{\Delta_\omega :
\omega\in \G\cap \supp \PP\}$, where $\Delta_\omega$ is the drift from a site with environment~$\omega$.
Suppose also that there is ${\hat \eps}_0$ such that
\[
 \PP\big[{\mathtt P}_{\!\text{\boldmath ${\omega}$}}[\text{total number of particles at time 1 is 1}]
               \geq {\hat \eps}_0\big]=1,
\]
i.e., for almost all environments, in any site the particle does not branch with uniformly 
positive probability. Then
\begin{equation}
\label{brwre_low}
 {\mathtt P}_{\!\text{\boldmath ${\omega}$}}[T(0,x_0)>n] \geq 
     \exp\big(-{\widehat K}_3^q(d) n \exp(-{\widehat K}^q_4(d) \ln^{1/d} n)\big).
\end{equation}
\end{prop}

Now, we go back to the random walk among soft obstacles.
In the one-dimensional case, we are able to obtain finer results. We assume now that the transition probabilities and the obstacles are independent, in other words,  $\PP=\mu\otimes\nu$, where
$\mu,\nu$ are two product measures governing, respectively, the transition
probabilities and the obstacles. We further assume ${\widetilde m}= 1$
and we denote $\omega_x^+ = \omega_x(+1)$ and $\omega_x^- =  1
-\omega_x^+ = \omega_x(-1)$.
Condition~N is now
 equivalent to
\begin{equation}
\label{Nbeta}
\begin{array}{r}
\inf\{a: \; \mu(\omega_0^ +   \le a)>0\}<1/2,\\
\sup\{a: \; \mu(\omega_0^+    \ge a)>0\}>1/2,
\end{array}
\end{equation}
i.e., the RWRE is strictly nestling.
Let
\begin{equation}
\label{rhodef}
\rho_i = \frac{\omega_i^-}{\omega_i^+},\quad i \in \Z\, .
\end{equation}
Define
$\kappa_\ell = \kappa_\ell(p)$ such that
\begin{equation}
\label{kappaldef}
\EE \Big(\frac{1}{\rho_0^{\kappa_\ell}}\Big) = \frac{1}{1-p}
\end{equation}
and 
$\kappa_r = \kappa _r(p)$ such that
\begin{equation}
\label{kappardef}
\EE (\rho_0^{\kappa_r}) = \frac{1}{1-p}
\end{equation}
Due to Condition~N, since $0 < p < 1$, $\kappa_\ell$ and $\kappa_r$ are well-defined, 
strictly positive and finite. Indeed, to see this for~$\kappa_r$,
observe that for the function $f(x)=\EE(\rho_0^x)$
it holds that $f(0)=1$, $f(x)\to\infty$ as $x\to\infty$, and~$f$ is convex, so the equation
$f(x)=u$ has a unique solution for any $u>1$. A similar 
argument implies that~$\kappa_\ell$
is well-defined.

We now  are able to characterize the
quenched and annealed tails of~$\tau$ in the following way:
\begin{theo}
\label{t_dim1_ann}
For $d=1$
 \begin{equation}
\label{dim1_ann}
 \lim_{n\to\infty}\frac{\ln\PA[\tau>n]}{\ln n} = -(\kappa_\ell (p) + \kappa_r (p)).
\end{equation}
\end{theo}

\begin{theo}
\label{dim_1_qu}
For $d=1$
\begin{equation}
\label{qu}
 \lim_{n\to\infty}\frac{\ln(-\ln \Ps[\tau>n])}{\ln n} = \frac{\kappa_\ell (p) + \kappa_r (p)}
{1+ \kappa_\ell (p) + \kappa_r (p)}
   \qquad\text{ $\PP$-a.s.}
\end{equation}
\end{theo}

In our situation, besides the quenched and the annealed
probabilities, one can also consider two ``mixed'' ones: the probability measure $\Po^z = \nu \times \Ps^z$ which is quenched in~$\om$ and
annealed in $\tb$ and the probability measure $\Pt^z = \mu\times \Ps^z$ which is quenched in~$\tb$ and
annealed in~$\om$.
Again, we use the simplified notations $\Po=\Po^0$, $\Pt=\Pt^0$.

Let 
\begin{align*}
 \beta_0 &=\inf\{\eps: \; \mu[\omega_0^+<\eps]>0\}\\
\beta_1 &=1-\sup\{\eps: \; \mu[\omega_0^+>\eps]>0\}.
\end{align*}
Due to \eqref{Nbeta} we have $ \beta_0 <1/2$,  $\beta_1 < 1/2 $.
Then, we have the following results about the ``mixed'' probabilities of survival:
\begin{theo}
\label{fixed_obstacles}
For $d=1$,
\begin{equation}
\label{qu_obst}
 \lim_{n\to\infty}\frac{\ln(-\ln \Pt[\tau>n])}{\ln n} = \frac{|\ln(1-p)|}{|\ln(1-p)|+ F_e},
   \qquad\text{ $\nu$-a.s.,}
\end{equation}
where 
\begin{equation}\label{Fedef}
F_e=\frac{\ln\Big(\frac{1-\beta_1}{\beta_1}\Big)\ln\Big(\frac{1-\beta_0}{\beta_0}\Big)}{\ln\Big(\frac{1-\beta_1}{\beta_1}\Big)+\ln\Big(\frac{1-\beta_0}{\beta_0}\Big)}.
\end{equation}
\end{theo}

\begin{theo}
 \label{fixed_omega}
For $d=1$, we have: 
\begin{itemize}
 \item[(i)] If $\EE(\ln\rho_0 )=0$, then, for each $\eps > 0$, there exist sequences of positive random variables $R_n^\eps(\om)$, $R_n(\om)$ and constants $K_1, K_2$ such that for $\mu$-almost all~$\om$,
\begin{align}
 \label{for_sinai_RW}
\Po[\tau>n]&\le e^{-K_1 R_n^\eps(\om) },\\
\Po[\tau>n]&\ge e^{-K_2 R_n(\om) }.
\end{align}
These random variables have the following properties:
 there exists a family of nondegenerate
random variables $(\Xi^{(\eps)},\eps\geq 0)$ such that
\[
\frac{R_n(\om)}{\ln^2 n}\to \Xi^{(0)} \qquad \text{and} \qquad 
\frac{R_n^\eps(\om)}{\ln^2n} \to \Xi^{(\eps)}
\]
in law as $n\to \infty$. Also, we have $\Xi^{(\eps)} \to \Xi^{(0)}$ in law as $\eps\to 0$. 
\item[(ii)]  If $\EE(\ln\rho_0)\ne 0$, then
\begin{equation}
 \label{not_sinai} \lim_{n\to\infty}\frac{\ln(-\ln \Po[\tau>n])}{\ln n} =\frac{\kappa}{\kappa+1},
\end{equation}
where $\kappa$ is such that $\EE(\rho_0^\kappa)=1$, in the case when $\EE(\ln\rho_0)< 0$, 
or $\EE\big(\frac{1}{\rho_0^\kappa}\big)=1$, in the case when $\EE(\ln\rho_0)>0$.
\end{itemize}
\end{theo}

\begin{rmk}
\label{rem_HuShi}
In fact, a comparison with Theorems~1.3 and~1.4 of~\cite{HS} suggests that
\begin{equation}
\label{asa}
\limsup\limits_{ n\to \infty} \frac{R_n^\eps(\om)}{\ln^2n \ln\ln \ln n} < \infty, \quad \liminf\limits_{ n\to \infty}\frac{R_n(\om)\ln \ln \ln n }{\ln^2n} > 0,
\end{equation}
and, in particular, for $\mu$-almost all $\omega$ and some positive constants $C_1,C_2$,
\begin{align} 
\label{is_sinai} 
&\limsup_{n\to\infty}\frac{\ln \Po[\tau>n]}{\ln^2 n \ln \ln \ln n} \leq -C_1, \\
\label{is_sinai2} 
&\liminf_{n\to\infty}\frac{\ln \Po[\tau>n]\ln \ln \ln n }{\ln^2 n }  \geq -C_2.
\end{align}
However, the proof of~(\ref{asa})--(\ref{is_sinai2}) would require a lengthy analysis
of fine properties of the potential~$V$ (see Definition~\ref{df_V} below), so we
decided that it would be unnatural to include it in this paper.
\end{rmk}
\begin{rmk}
It is interesting to note that~$r$ does only enter the constants, but not the exponents in all these results.
\end{rmk}
\begin{rmk}
\label{moge}
In fact, the proofs of~(\ref{up_ann}) and~(\ref{up_qu_d}) do not really use independence (and can be readily extended
to the finitely dependent case), 
but we will use independence for the proofs of the lower bounds in~(\ref{prob_trap}).
However, if one modifies Condition~N in a suitable way, we conjecture that Theorem~\ref{t_ann} and Theorem~\ref{t_qu} remain true if the environment is not i.i.d.\ but finitely dependent (for BRWRE, one can find generalizations of this kind in~\cite{CP07-2}).
\end{rmk}

\section{Proofs: multi-dimensional case}
\label{s_proofs_d}
In this section, we prove Theorems~\ref{t_ann} and~\ref{t_qu}.
In fact,            
the ideas we need to prove these results are
similar to those in the  proofs of Theorems~1.8 and~1.9 of~\cite{CP07-1}.
In the following, we explain the relationship of the discussion in~\cite{CP07-1}
with our model, and give the proof of Theorems~\ref{t_ann} and~\ref{t_qu},
sometimes referring to~\cite{CP07-1} for a more detailed account.

\medskip
\noindent
{\it Proof of Theorems~\ref{t_ann} and~\ref{t_qu}.}
The proof of~(\ref{up_ann}) follows essentially the proof of Theorem~1.8
of~\cite{CP07-1}, where it is shown that the tail 
of the first hitting time of some fixed site~$x_0$
(one may think also of the first return time to the origin) can be bounded
 from above as in~(\ref{up_ann}). 
The main idea is that, as a general fact, for any recurrent BRWRE there are the 
so-called \emph{recurrent seeds}. These are simply finite configurations
of the environment, where, with positive probability, the number of particles
grows exponentially without help from outside (i.e., suppose that all
particles that step outside this finite piece are killed; then, the number
of particles in the seed dominates a supercritical Galton-Watson process, which
explodes with positive probability). Then, we consider an embedded RWRE,
until it hits a recurrent seed and the supercritical Galton-Watson process
there explodes (afterwards, the particles created by this explosion are used
to find the site~$x_0$, but here this part is only needed for Proposition \ref{brwre}).

So, going back to the model of this paper, 
obstacles play the role of 
recurrent seeds, and the moment~$\tau$ when the event $\{\xi_n\in\Theta,Z_n=1\}$ happens for the first time
is analogous to the moment of the first explosion of the Galton-Watson process in the recurrent seed. To explain better this analogy, consider the following
situation. Suppose that, outside the recurrent seeds there is typically a strong drift 
in one direction and the branching is very weak or absent. Then, the qualitative behaviour
of the process is quite different before and after the first explosion. 
Before, we typically observe very few (possibly even one) particles with more or less
ballistic behaviour; after, the cloud of particles starts to grow exponentially
in one place (in the recurrent seed where the explosion occurs), and so the cloud of particles
expands linearly in all directions. So, the first explosion of one of the Galton-Watson 
processes in recurrent seeds marks the transition between qualitatively different behaviours
of the BRWRE, and thus it is analogous to the moment~$\tau$ of the model of
the present paper.

First, we prove~(\ref{up_ann}) for $d\geq 2$. For any $a\in \Z$, define $\KK_a=[-a, a]^d$.
Choose any $\alpha<(\ln \eps_0^{-1})^{-1}$ ($\eps_0$ is from Condition~E) and define the event
\[
 M_n = \{\s : \text{ for any }y\in\KK_{{\widetilde m}n} \text{ there exists }
 z\in\Theta \text{ such that }\|y-z\|_1\leq \alpha\ln n\}
\]
(recall that ${\widetilde m}$ is a constant such that $\omega_x(y)=0$ if
$\|x-y\|>{\widetilde m}$, introduced in Section~\ref{s_intro}).
Clearly, we have for $d\geq 2$
\begin{equation}
\label{oc_Mn_compl}
 \PP[M_n^c] \leq C_1n^d\exp(-C_2 \ln^d n).
\end{equation}

Now, suppose that $\s\in M_n$. So, for any possible location of the particle 
up to time~$n$, we can find a site with an obstacle which is not more than
$\alpha\ln n$ away from that location (in the sense of $L^1$-distance).
 This means that, on any time interval
of length $\alpha\ln n$, the particle will disappear (i.e., $\tau$ is in this interval if the particle has not disappeared before)
with probability at least $r\eps_0^{\alpha\ln n}$, where $\eps_0$
is the constant from the uniform ellipticity condition.
There are $\frac{n}{\alpha\ln n}$ such (disjoint) intervals in the interval $[0,n]$, so 
\begin{align}
 \Ps[\tau>n] &\leq (1-r\eps_0^{\alpha\ln n})^{\frac{n}{\alpha\ln n}} \nonumber\\
 &\leq \exp\Big(-\frac{C_3 n^{1-\alpha\ln\eps_0^{-1}}}{\ln n}\Big).\label{up_qu_rough}
\end{align}
Then, from~(\ref{oc_Mn_compl}) and~(\ref{up_qu_rough}) we obtain 
(recall that $\alpha\ln\eps_0^{-1}<1$)
\[
 \PA[\tau>n] \leq \exp\Big(-\frac{C_3 n^{1-\alpha\ln\eps_0^{-1}}}{\ln n}\Big)
 + C_1n^d\exp(-C_2 \ln^d n),
\]
and hence~(\ref{up_ann}).

Let us now prove~(\ref{up_qu_d}), again for the case $d\geq 2$.
 Abbreviate by $\ell_d=\frac{2^d}{d!}$ the volume
of the unit sphere in~$\R^d$ with respect to the $L^1$ norm, and let $q=\PP[\theta_0=0]= 1-p$.
Choose a large enough~$\widehat \alpha$ in such a way that $\ell_d {\widehat\alpha}^d\ln q^{-1}>d+1$,
and define
\begin{align*}
 {\widehat M}_n &= \{\s : \text{ for any }y\in\KK_{{\widetilde m}n} \text{ there exists }
 z\in\Theta \\
  & \qquad\qquad\qquad\text{ such that }\|y-z\|_1\leq \widehat\alpha\ln^{1/d} n\}.
\end{align*}
By a straightforward calculation, we obtain for $d\geq 2$
\begin{equation}
\label{oc_hat_M}
 \PP[{\widehat M}_n^c] \leq C_4n^d n^{-\ell_d {\widehat \alpha}^d\ln q^{-1}}.
\end{equation}
 Using the Borel-Cantelli lemma, (\ref{oc_hat_M}) implies that for $\PP$-almost all~$\s$,
there exists~$n_0(\s)$ such that $\s\in {\widehat M}_n$ for all $n\geq n_0(\s)$.

Consider now an environment~$\s \in {\widehat M}_n$.
In such an environment, in the $L^1$-sphere of size~$n$ around the origin, 
any $L^1$-ball of radius~$\widehat{\alpha} \ln^{1/d}n$ contains at least 
one obstacle (i.e., a point from~$\Theta$). This means that, in any time interval
of length $\widehat\alpha\ln^{1/d}n$, the particle will disappear
with probability at least $r\eps_0^{\widehat\alpha\ln^{1/d}n}$, where, as before, $\eps_0$
is the constant from the uniform ellipticity Condition~E.
There are $\frac{n}{\widehat\alpha\ln^{1/d}n}$ such intervals on $[0,n]$, so
\[
 \Ps[\tau>n] \leq (1-r\eps_0^{\widehat\alpha\ln^{1/d}n})^{\frac{n}{\widehat\alpha\ln^{1/d}n}},
\]
which gives us~(\ref{up_qu_d}) in dimension~$d\geq 2$.

Now, we obtain~(\ref{up_ann}) and~(\ref{up_qu_d})
in the one-dimensional case. Since the environment is i.i.d., there
 exist~$\gamma_1,\gamma_2>0$ such that for any interval~$I\subset\Z$,
\begin{equation}
\label{mnogo_tochek}
 \PP[|I\cap\Theta| \geq \gamma_1|I|] \geq 1-e^{-\gamma_2|I|}.
\end{equation}
We say that an interval~$I$ is \emph{nice}, if it
contains at least~$\gamma_1 |I|$ sites from~$\Theta$.

Define
\begin{align*}
 h(\s) &= \min\{\text{$m$: all the intervals of length $k\geq m$} \\
 & ~~~~~~~~~~~~~~~~~~~~~~
 \text{intersecting with $[-e^{\gamma_2k/2},e^{\gamma_2k/2}]$ are nice}\}.
\end{align*}
It is straightforward to obtain from~(\ref{mnogo_tochek}) that
there exists~$C_5>0$ such that
\begin{equation}
\label{hvost_odnako}
 \PP[h(\s)>  k] \leq e^{-C_5 k}.
\end{equation}
In particular, $h(\s)$ is finite for $\PP$-almost all~$\s$.

Now, define the event $F=\{\max_{s\leq n}|\xi_s| \leq n^{a}\}$,
where $a = \frac{\gamma_2}{4\ln \eps_0^{-1}}$. By Condition~E, 
during any time interval of length $\frac{\ln n}{2\ln \eps_0^{-1}}$ the random walk 
completely covers a space interval of the same length with probability at
least
\[
 \eps_0^{\frac{\ln n}{2\ln \eps_0^{-1}}} = n^{-1/2}.
\]
Assume that $h(\s) < \frac{\ln n}{2\ln \eps_0^{-1}}$ and
observe that $\frac{\gamma_2}{2}\times \frac{1}{2\ln \eps_0^{-1}} = a$,
so that all the intervals of length at least $\frac{\ln n}{2\ln \eps_0^{-1}}$
intersecting with $[-n^a,n^a]$ are nice. On the event $F$,
 the random walk visits the set~$\Theta$ at least $O(\frac{n^{1/2}}{\ln n})$
times with probability at least $1-\exp(-C_6\frac{n^{1/2}}{\ln n})$. Indeed, split the time into $\frac{2 n \ln\eps_0^{-1}}{\ln n}$ intervals of length $\frac{\ln n}{2\ln \eps_0^{-1}}$, and consider such an interval successful if the random walk 
completely covers a space interval of the same length: we then have $\frac{2 n \ln\eps_0^{-1}}{\ln n}$ independent trials with success probability 
at least~$n^{-1/2}$, and then one can use Chernoff's bound for Binomial distribution (see e.g.\ inequality~(34) of~\cite{CP07-1}).
Hence we obtain for such~$\s$
\begin{equation}
\label{trapped_blizko}
 \Ps[\tau>n, F] \leq \exp\Big(-C_7\frac{n^{1/2}}{\ln n}\Big).
\end{equation}

Let us define the sequence of stopping times $t_k$, $k=0,1,2,\ldots$
as follows: $t_0=0$ and
\[
 t_{k+1} = \min\{t\geq t_k+{\widetilde m} : \text{ there exists }
   z\in\Theta\text{ such that }|\xi_t-z|\leq {\widetilde m}-1\}
\]
for $k\geq 1$. Defining also the sequence
of events
\[
 D_k = \{\text{there exists }t\in[t_k,t_k+{\widetilde m}-1] \text{ such that }
               \xi_t\in\Theta \text{ and }Z_t=1\},
\]
by Condition~E we have 
\begin{equation}
\label{popalso}
 \Ps[D_{k+1}\mid \FFF_k] \geq r\eps_0^{{\widetilde m}-1},
\end{equation}
where~$\FFF_k$ is the sigma-algebra generated by 
$D_0,\ldots,D_k$. 

Observe that, for~$\s$ with $h(\s) < \frac{\ln n}{2\ln \eps_0^{-1}}$,
on the event~$F^c$ we have $t_{k'}<n$, where $k'=\frac{\gamma_1}{2{\widetilde m}^2}n^{a}$.
Thus, by~(\ref{popalso}),
\begin{equation}
\label{trapped_daleko}
 \Ps[\tau>n, F^c] \leq \exp(-C_8 n^{a}),
\end{equation}
and we obtain~(\ref{up_qu_d}) from~(\ref{trapped_blizko}) and~(\ref{trapped_daleko}) (notice that in the one-dimensional
case, the right-hand side of~\eqref{up_qu_d} is of the form
$\exp(-K_1^q(1)n^{1-K_2^q(1)})$).
Then, the annealed upper bound~(\ref{up_ann}) for $d=1$ follows from~(\ref{up_qu_d})
and~(\ref{hvost_odnako}).

Now, let us prove the lower bound~(\ref{low_ann}).
This time, we proceed as in the proof of Theorem~1.9 from~\cite{CP07-1}.
Denote by $\Sph^{d-1} = \{x\in\R^d : \|x\|_2=1\}$ the unit sphere in~$\R^d$,
 and, recalling (\ref{Deltadef}), let~$\Delta_\omega$ be the drift 
at the point $(\omega,\theta)\in\MM$.
One can split the sphere $\Sph^{d-1}$ into a finite number
(say, $m_0$) of non-inter\-sec\-ting subsets $U_1,\ldots, U_{m_0}$
and find a finite collection 
$\Gamma_1,\ldots, \Gamma_{m_0}\subseteq \MM$
having the following properties: for all $i=1,\ldots,m_0$,
\begin{itemize}
\item[(i)] $\theta=0$ for all $\sigma=(\omega,\theta)\in \Gamma_i$,
\item[(ii)] there exists~$p_1>0$ (depending only on the law of the environment)
 such that $\PP[\sigma_0\in\Gamma_i]>p_1$,
\item[(iii)] there exists~$a_1>0$ such that
for any $z\in U_i$ and any $\sigma=(\omega,\theta)\in {\Gamma}_i$
we have $z\cdot \Delta_\omega < -a_1$ (recall Condition~N).
\end{itemize}
Intuitively, this collection will be used to construct (large)
pieces of the environment which are free of obstacles (item (i))
and have the drift pointing towards the center of the corresponding region
(item (iii)). The cost of constructing piece of environment of size~$N$
(i.e., containing~$N$ sites) with such properties does not exceed
$p_1^N$ (item (ii)).

Consider any $z\in\Z^d$, 
$B\subset \Z^d$ and a collection $H=(H_x\subseteq \MM,x\in B)$;
let us define
\[
 \SSS(z,B,H) = \{\s : \sigma_{z+x}\in H_x \text{ for all }x\in B\}.
\]
In~\cite{CP07-1}, on $\SSS(z,B,H)$ we said that there is an $(B,H)$-seed
in~$z$; for the model of this paper, however, we prefer not to use
the term ``seed'', since the role seeds play in~\cite{CP07-1} is
quite different from the use of environments belonging to $\SSS(z,B,H)$ here.
Take $G^{(n)}=\{y\in\Z^d: \|y\|_2 \leq u\ln n\}$, 
where~$u$ is a (large) constant to be chosen later. 
Let us define the sets $H^{(n)}_x,x\in G^{(n)}$ in the following way.
First, put $H^{(n)}_0=\Gamma_1$; for $x\neq 0$, let~$i_0$ be such that
$\frac{x}{\|x\|}\in U_{i_0}$ (note that~$i_0$ is uniquely
defined), then put $H^{(n)}_x=\Gamma_{i_0}$.
Clearly, for any $y\in\Z^d$
\begin{equation}
\label{prob_trap}
 \PP[\SSS(y,G^{(n)},H^{(n)})] \geq p_1^{(2u)^d \ln^d n}.
\end{equation}

Denote
\[
 {\widehat p} = \sup_{y:\|y\|_1\leq {\widetilde m}} \Ps^y[\xi\text{ hits }\Z^d\setminus G^{(n)}
  \text{ before }0].
\]
As in~\cite{CP07-1} (see the derivation of~(42) there), we obtain that there
exist $C_9,C_{10}$ such that for all $\s\in\SSS(0,G^{(n)},H^{(n)})$ we have
\begin{equation}
\label{prob_escape_d}
 {\widehat p} \leq \frac{C_9}{n^{C_{10}u}}.
\end{equation}
So, choose $u>\frac{1}{C_{10}}$, then, 
on the event that~$\s\in\SSS(0,G^{(n)},H^{(n)})$, (\ref{prob_escape_d}) implies that
\begin{equation}
\label{staylong}
 \Ps[\xi_i\in G^{(n)} \text{ for all }i\leq n] \geq (1-{\widehat p})^n \geq C_{11}
\end{equation}
(if the random walk hits the origin at least~$n$ times before
hitting $\Z^d\setminus G^{(n)}$, then $\xi_i\in G^{(n)}$ 
for all $i\leq n$).
Since~$G^{(n)}$ is free of obstacles, we obtain~(\ref{low_ann}) from~(\ref{prob_trap}) and (\ref{staylong}).

Now, it remains to prove~(\ref{low_qu_d}). Define
\[
 {\widehat G}^{(n)} = \{y\in\Z^d : \|y\|_2 \leq v \ln^{1/d}n\},
\]
and let ${\widehat H}^{(n)} = ({\widehat H}^{(n)}_x,x\in {\widehat G}^{(n)})$ be defined in the same way as~$H^{(n)}$ above,
but with~${\widehat G}^{(n)}$ instead of~$G^{(n)}$. Analogously to~(\ref{prob_trap}), we have
\begin{equation}
\label{prob_trap2}
 \PP[\SSS(y,{\widehat G}^{(n)},{\widehat H}^{(n)})] \geq p_1^{(2v)^d \ln n}
                          = n^{-(2v)^d \ln p_1^{-1}}.
\end{equation}
Choose~$v$ in such a way that $b_0:=(2v)^d \ln p_1^{-1}<d/2(d+1)$.
Then, it is not difficult to obtain (by dividing $\KK_{\sqrt{n}}$ into
$O(n^{d(\frac{1}{2}-b_0)})$ subcubes of linear size~$O(n^{b_0})$) that
\begin{equation}
\label{nashli_lovushko_d}
 \PP\Big[\bigcup_{z\in\KK_{\sqrt{n}}}\SSS(z,{\widehat G}^{(n)},
              {\widehat H}^{(n)})\Big] \geq
    1-\exp(-C_{12}n^{\frac{d}{2}-(d+1) b_0)}).
\end{equation}
Using the Borel-Cantelli Lemma, $\PP$-a.s.\ for all~$n$ large enough 
we have 
\[
\s\in \bigcup_{z\in\KK_{\sqrt{n}}}\SSS(z,{\widehat G}^{(n)},{\widehat H}^{(n)}).
\]
Denote by~$T_B$ the first hitting time of a set~$B\subset\Z^d$:
\[
 T_B = \inf\{m\geq 1: \xi_m\in B\},
\]
and write $T_a=T_{\{a\}}$ for one-point sets. Next,
for $\s\in \SSS(0,{\widehat G}^{(n)},{\widehat H}^{(n)})$ we are going
to obtain an upper bound for 
$q_x:= \Ps^x[\xi\text{ hits }\Z^d\setminus {\widehat G}^{(n)}
  \text{ before }0]
=
\Ps^x[T_0>T_{\Z^d\setminus {\widehat G}^{(n)}}]$, uniformly
in~$x\in {\widehat G}^{(n)}$. To do this, note that there are positive constants $C_{13},C_{14}$
such that, abbreviating $B_0=\{x\in\Z^d: \|x\|_2\leq C_{13}\}$, the process
$\exp(C_{14}\|\xi_{m\wedge T_{B_0}}\|_2)$ is a supermartingale (cf.\ the proof
of Theorem~1.9 in~\cite{CP07-1}), i.e.,
\[
 \Es\big(\exp(C_{14}\|\xi_{(m+1)\wedge T_{B_0}}\|_2)\mid 
\xi_{j\wedge T_{B_0}}, j\leq m\big)
 \leq \exp(C_{14}\|\xi_{m\wedge T_{B_0}}\|_2)
\]
for all~$m$.
 Denote
\[
 {\widetilde G}^{(n)} = \{x\in\Z^d : \|x\|_2 \leq v \ln^{1/d}n - 1\}.
\]
For any $x\in {\widetilde G}^{(n)}$ and $y\in \Z^d\setminus {\widehat G}^{(n)}$,  we have
$\|x\|_2\leq \|y\|_2-1$, so $e^{C_{14}\|x\|_2}\leq e^{-C_{14}}e^{C_{14}\|y\|_2}$.
Keeping this in mind, we apply the Optional Stopping Theorem to obtain that,
for any $\s\in \SSS(0,{\widehat G}^{(n)},{\widehat H}^{(n)})$,
\[
 q_x \exp(C_{14} v \ln^{1/d}n) \leq \exp(C_{14} \|x\|_2) \leq e^{-C_{14}} \exp(C_{14} v \ln^{1/d}n),
\]
so $q_x\leq e^{-C_{14}}$ for all $x \in \widetilde G$. Now, from any $y\in {\widehat G}^{(n)}\setminus {\widetilde G}^{(n)}$ the particle 
can come to~${\widetilde G}^{(n)}$ in a fixed number of steps (at most $\sqrt{d}+1$) with 
uniformly positive probability. This means that,
on $\SSS(0,{\widehat G}^{(n)},{\widehat H}^{(n)})$,
 there exists a positive constant $C_{15}>0$ such that
for all $x\in {\widehat G}^{(n)}$
\begin{equation}
\label{escapar-nicht!}
 \Ps^x[T_0<T_{\Z^d\setminus {\widehat G}^{(n)}}] \geq C_{15}.
\end{equation}
Then, analogously to~(\ref{prob_escape_d}), 
on $\SSS(0,{\widehat G}^{(n)},{\widehat H}^{(n)})$
we obtain that, for all~$y$ such that
$\|y\|_1\leq {\widetilde m}$
\[
 \Ps^y[\xi\text{ hits }\Z^d\setminus {\widehat G}^{(n)}
  \text{ before hitting }0] \leq \frac{C_{16}}{\ln^{1/d}n}.
\]
So, using~(\ref{escapar-nicht!}) on $\SSS(0,{\widehat G}^{(n)},{\widehat H}^{(n)})$
we obtain that there are~$C_{17}$ and~$C_{18}$ such that for all $x\in {\widehat G}^{(n)}$
\begin{equation}
\label{surv_log1/d}
 \Ps^x[T_{\Z^d\setminus {\widehat G}^{(n)}}\geq \exp(C_{17}\ln^{1/d}n)]\geq C_{18}.
\end{equation}
\begin{figure}
\centering
\includegraphics{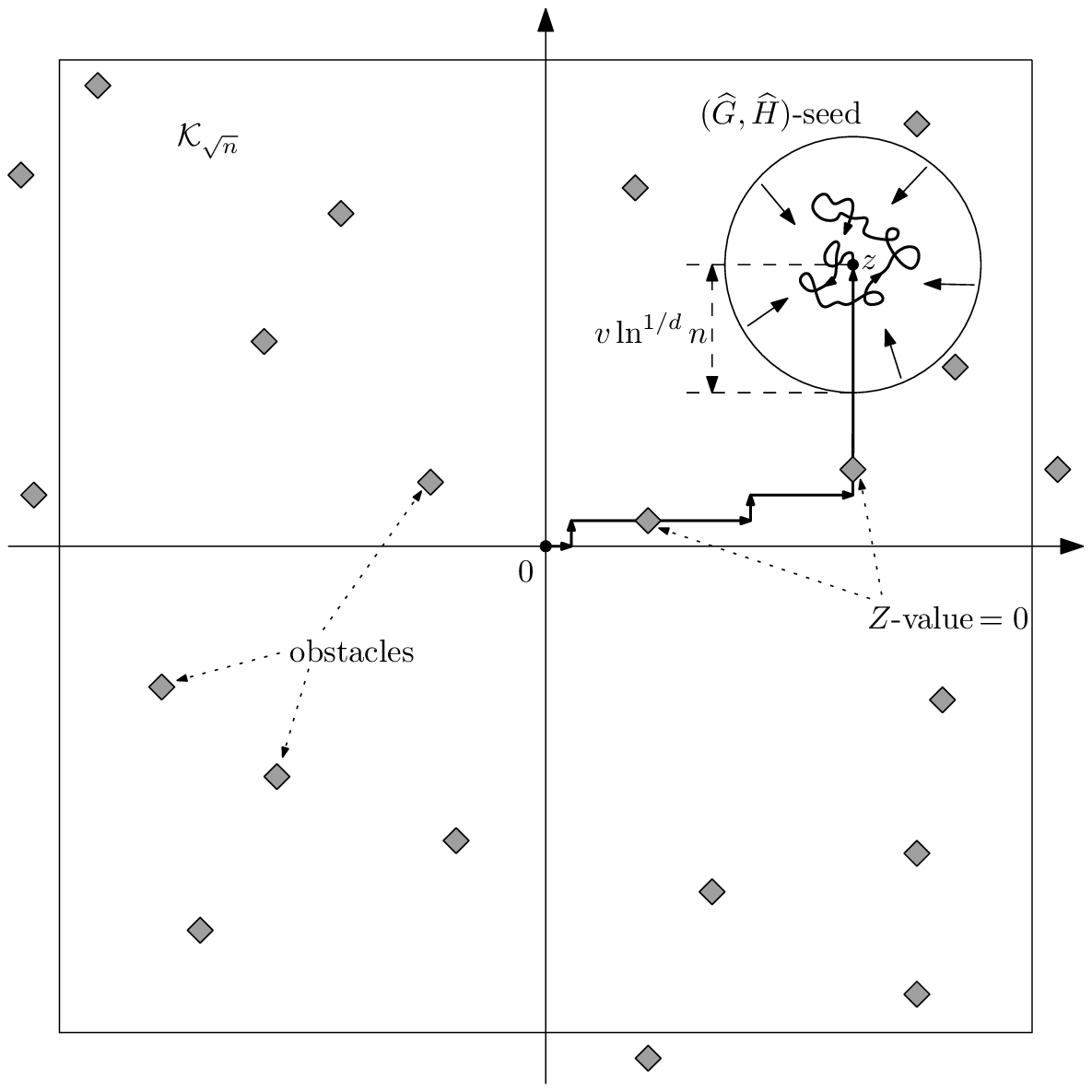}
\caption{The (quenched) strategy of survival used in the proof
of~(\ref{low_qu_d})}
\label{f_q_surv}
\end{figure}
Then, we use the following survival strategy of the particle (see Figure~\ref{f_q_surv}): 
provided that the event $\SSS(z,{\widehat G}^{(n)},{\widehat H}^{(n)})$
occurs for some $z\in\KK_{\sqrt{n}}$,
first the particle walks there (using, for instance, the shortest possible path)
without disappearing in an obstacle; this happens 
with probability at least $(\eps_0(1-r))^{d\sqrt{n}}$.
Then, it will spend time
\[
 n = \exp(C_{17}\ln^{1/d}n) \times n \exp(-C_{17}\ln^{1/d}n)
\]
in $z+{\widehat G}^{(n)}$ 
with probability at least $C_{18}^{n \exp(-C_{17}\ln^{1/d}n)}$, so
\[
 \Ps[\tau>n] \geq (\eps_0(1-r))^{d\sqrt{n}} C_{18}^{n \exp(-C_{17}\ln^{1/d}n)},
\]
and this gives us~(\ref{low_qu_d}).
\qed

\medskip
\noindent
{\it Proof of Proposition~\ref{brwre}.}
Now, we explain how to obtain Proposition~\ref{brwre}. To prove~\eqref{brwre_up}, we
proceed as in the proof of~\eqref{up_qu_d}. As noted in the beginning of this section,
the disappearing of the particle in an obstacle is analogous to starting an exploding Galton-Watson
process in a recurrent seed. 
Denote by~$\widetilde T$ the moment when this happens, i.e., at least
$e^{C_{19}k}$ new particles are created in this recurrent seed by time
${\widetilde T}+k$. Thus, one can obtain a bound of the form 
\[
 {\mathtt P}_{\text{\boldmath ${\omega}$}}[{\widetilde T}>n] \leq \exp\big(-C_{20} n \exp(-C_{21} \ln^{1/d} n)\big).
\]
Then, using the uniform ellipticity, it is straightforward to obtain that, waiting $C_{22}n$
time units more (with large enough $C_{22}$), one of the newly created (in this recurrent seed)
particles will hit~$x_0$ with probability at least $1-e^{-C_{23}n}$, and this implies~\eqref{brwre_up}.

To show~\eqref{brwre_low}, we note that, analogously to the proof of~(\ref{low_qu_d}) that we are able
to create a seed which is free of branching sites of diameter $C_{24}\ln^{1/d}n$, which lies
at distance $O(\sqrt{n})$ from the origin. Then, the same idea works: the initial particle goes
straight to the seed without creating new particles, and then stays there up to time~$n$.
The detailed proof goes along the lines of the proof of~(\ref{low_qu_d}) only with 
notational adjustments.
\qed

\section{Proofs: one-dimensional case}
\label{s_proofs_1}
\subsection{Preliminaries} 
We define the potential, which is a function of the transition probabilities.
 Under our assumptions it is a sum of i.i.d.\ random variables. Recall~(\ref{rhodef}).
\begin{df}
\label{df_V}
Given the realization of the random environment,
the potential~$V$ is defined by
\[
V(x)=\left\{\begin{array}{ll}
\sum_{i=1}^{x}\ln\rho_i, & x>0,\\
0, & x=0,\\
\sum_{i=x+1}^{0}\ln\frac{1}{\rho_i}, & x<0.
\end{array}
 \right.
\]
 \end{df}

\begin{df}
\label{trap}
We say that there is a trap of depth $h$ located at 
$ [x-b_1, x+b_2]$ with the bottom at $x$ if
\begin{align*}
&V(x)=\min_{y\in [x-b_1, x+b_2]} V(y)\\
&V( x-b_1)-V(x)\ge h\\
&V( x+b_2)-V(x)\ge h\, .
\end{align*}
Note that we actually require the depth of the trap to be \emph{at least} $h$.  
We say that the trap is \emph{free of obstacles} 
if in addition $\Theta\cap [x-b_1, x+b_2]=\emptyset$. 
\end{df}
Define
\begin{equation}
\label{psidef}
\psi(h,b_1, b_2)=
 \sup_{\lambda>0}\Big\{\lambda h-b_2\ln \EE (\rho_0^\lambda)\Big\}
+\sup_{\lambda>0}\Big\{\lambda h-
     b_1\ln \EE \Big(\frac{1}{\rho_0^\lambda}\Big)\Big\}
 \end{equation}
and 
\begin{equation}
\label{psitilddef}
\widetilde \psi(h)=\inf_{b_1, b_2>0}\{-(b_1+b_2)\ln (1-p)+\psi(h, b_1, b_2)\}.
\end{equation}

\begin{lmm} 
\label{traps}
Let $\Lambda^x(h, b_1, b_2, n)$ be  the event that there is a trap of depth $h\ln n$,
located  at $ [x-b_1\ln n, x+b_2\ln n]$ with the bottom at~$x$. Let also
$A^x(h, b_1, b_2, n)$ be the event that there is a trap of depth $h\ln n$, \emph{free of obstacles},
located  at $ [x-b_1\ln n, x+b_2\ln n]$ with the bottom at~$x$. For any $0<h^-, h^+, b_i^-, b_i^+<\infty$, we have 
\begin{align}
\label{eq_sem_traps}
\lim_{n\to\infty}\frac{\ln \PP[\Lambda^x(h, b_1, b_2, n)]}{\ln n}
&=-\psi(h, b_1, b_2),\\
\label{eq_traps}
\lim_{n\to\infty}\frac{\ln \PP[A^x(h, b_1, b_2, n)]}{\ln n}
&= (b_1+ b_2)\ln (1-p) - \psi(h, b_1, b_2),
\end{align}
uniformly in $h\in (h^-, h^+)$, $b_i\in (b_i^-, b_i^+)$.
\end{lmm}

\medskip
\noindent
{\it Proof of Lemma~\ref{traps}}. 
Note that the part $(b_1+ b_2)\ln (1-p)$ in \eqref{eq_traps} corresponds to the probability that 
the interval $ [x-b_1\ln n, x+b_2\ln n]$ is obstacle free. 

For notational convenience, we often omit integer parts and write, e.g., $b_1\ln n$ instead of its integer part. 
Using Chebyshev's inequality, we have for $\lambda > 0$
\begin{align*}
\PP[V( x+b_2\ln n)-V(x)\ge h \ln n]
&= \PP\Big[\exp\Big(\lambda\sum_{i=1}^{b_2\ln n}\ln \rho_i\Big)\ge \exp(\lambda h\ln n)\Big]\\
&\le \exp\big(-(\ln n)(\lambda h-b_2\ln [\EE e^{\lambda \ln \rho_0}])\big).
\end{align*}
Thus, 
\begin{align*}
\PP[V( x+b_2\ln n)-V(x)\ge h \ln n]
&\le \exp\Big(-(\ln n) \sup_{\lambda>0}\Big\{\lambda h-
     b_2\ln \EE (\rho_0^\lambda)\Big\}\Big).
\end{align*}
 
Analogously, we obtain
\begin{align*}
\PP[V( x-b_1\ln n)-V(x)\ge h \ln n]
&\le \exp\!\Big(-\!(\ln n) \sup_{\lambda>0}\!\Big\{\lambda h-
     b_1\ln \EE \Big(\frac{1}{\rho_0^\lambda}\Big) \Big\}\!\Big).
\end{align*}
So, 
\begin{align*}
\PP[A^x(h, b_1, b_2, n)]
&\le (1-p)^{(b_1+b_2)\ln n}\\
& \quad \times
 \exp\Big(-(\ln n) \Big[\sup_{\lambda>0}\Big\{\lambda h-
     b_2\ln \EE (\rho_0^\lambda)\Big\} \\
   &~~~~~~~~~~~~~~~~~~~~+ \sup_{\lambda>0}\Big\{\lambda h-
     b_1\ln \EE \Big(\frac{1}{\rho_0^\lambda}\Big) \Big\}\Big]\Big).
\end{align*}
and 
\begin{align*}
\PP[\Lambda^x(h, b_1, b_2, n)]
&\le 
 \exp\Big(-(\ln n) \Big[\sup_{\lambda>0}\Big\{\lambda h-
     b_2\ln \EE (\rho_0^\lambda)\Big\} \\
   &~~~~~~~~~~~~~~~~~~~~+ \sup_{\lambda>0}\Big\{\lambda h-
     b_1\ln \EE \Big(\frac{1}{\rho_0^\lambda}\Big) \Big\}\Big]\Big).
\end{align*}
To show~\eqref{eq_sem_traps} and~\eqref{eq_traps}, we have to obtain now the corresponding 
lower bounds. 
To this end, note first that, by 
Cram\'er's Theorem,
\begin{equation}
\label{LDP_rho}
 \lim_{k\to\infty}\frac{1}{b_2k}\ln \PP\Big[\sum_{i=1}^{b_2k}
\ln \rho_i\geq hk\Big]=
   \sup_{\lambda>0}\Big\{\lambda h - b_2\ln \EE(\rho_0^\lambda)\Big\}
\end{equation}
(recall that we treat $b_2k$ as an integer).

Define $S_\ell = \sum_{i=1}^\ell \ln \rho_i$, and, for $j\in [1,b_2k]$
\[
 S_\ell^{(j)} = \sum_{i=1}^\ell \ln \rho_{(i+j-2 \mod b_2 k)+1}.
\]
We have (recall that $h>0$)
\begin{align*}
\lefteqn{ b_2k\PP\Big[S_{b_2k}\geq hk, S_\ell\geq 0\text{ for all }
        \ell=1,\ldots,b_2k\Big]} \\
&= \sum_{j=1}^{b_2k} \PP\Big[S_{b_2k}\geq hk, 
       S_\ell^{(j)}\geq 0\text{ for all } \ell=1,\ldots,b_2k\Big]\\
 &\geq \PP\Big[S_{b_2k}\geq hk, \text{ there exists $j$ such that }
       S_\ell^{(j)}\geq 0\text{ for all } \ell=1,\ldots,b_2k\Big]\\
 &=\PP\Big[S_{b_2k}\geq hk\Big],
\end{align*}
since if $\sum_{i=1}^{b_2k}\ln \rho_i\geq hk$, then choosing~$j$ in such a
way that $S_j\leq S_\ell$ for all $\ell=1,\ldots,b_2k$, it is straightforward to
obtain that $S_\ell^{(j)}\geq 0$ for all $\ell=1,\ldots,b_2k$.

Hence
\[
 \PP\Big[S_{b_2k}\geq hk, S_\ell\geq 0\text{ for all }
        \ell=1,\ldots,k\Big] \geq \frac{1}{b_2k}\PP\Big[S_{b_2k}\geq hk\Big],
\]
which permits us to obtain a lower bound on
\[
 \PP[V(x+b_2\ln n)-V(x)>h\ln n, V(y)\geq V(x) \text{ for all }y\in(x,x+b_2\ln n)].
\]
Then, one obtains~\eqref{eq_sem_traps} and~\eqref{eq_traps} from~\eqref{LDP_rho} and the corresponding statement 
with~$b_1$ instead of~$b_2$ and~$1/\rho_i$ instead of~$\rho_i$.
\qed

\medskip

Next, we obtain a simpler expression for the function $\widetilde{\psi}$ (recall~\eqref{kappaldef}, \eqref{kappardef},
and~\eqref{psitilddef}).
\begin{lmm} 
\label{psitildesolution}
 We have
\begin{equation}
\label{psitilde}
\widetilde{\psi}(h) = (\kappa_\ell (p) + \kappa_r(p) ) h.
\end{equation}
\end{lmm}
\noindent
{\it Proof.}
By~(\ref{psidef}) and~(\ref{psitilddef}), it holds that
\begin{align*}
\widetilde{\psi}(h) =& 
\inf_{b_1 > 0}\Big\{ \sup_{\lambda>0}\Big\{\lambda h- b_1\ln \EE (1-p) \Big(\frac{1}{\rho_0^\lambda}\Big)\Big\}\Big\} \\
& {}+ \inf_{b_2 > 0}\Big\{ \sup_{\lambda>0}\Big\{\lambda h- b_2\ln \EE (1-p) (\rho_0^\lambda)\Big\}\Big\} \, .
\end{align*}
We will show that
\begin{equation}
\label{rightkap}
\inf_{b_2 > 0}\Big\{ \sup_{\lambda>0}\Big\{\lambda h- b_2\ln \EE (1-p) (\rho_0^\lambda)\Big\}\Big\} = \kappa_r h = \kappa_r(p)h.
\end{equation}
In the same way, one proves
\[
\inf_{b_1 > 0}\Big\{ \sup_{\lambda>0}\Big\{\lambda h- b_1\ln \EE (1-p) \Big(\frac{1}{\rho_0^\lambda}\Big)\Big\}\Big\} = \kappa_\ell(p)h\, .
\]
To show (\ref{rightkap}), note that
taking $\lambda = \kappa_r$ yields
\[
\inf_{b_2 > 0}\Big\{ \sup_{\lambda>0}\Big\{\lambda h- b_2\ln \EE (1-p)(\rho_0^\lambda)\Big\}\Big\} \geq \kappa_r h.
\]
Consider the function $g_b(\lambda) = \lambda h  - b\ln \EE \big((1-p)\rho_0^\lambda\big)$. Clearly, $g_b(0) = b\ln\frac{1}{1-p} > 0$
and an elementary calculation shows that for all $b\in (0, \infty)$
the function $g_b$ is concave (indeed, by the Cauchy-Schwarz inequality,
for any positive $\lambda_1,\lambda_2$ we obtain $\ln \EE\Bigl(\rho_0^{\frac{\lambda_1+\lambda_2}{2}}\Bigr) \leq \frac{1}{2}(\ln \EE(\rho_0^{\lambda_1})+
\ln \EE(\rho_0^{\lambda_2}))$).
Taking 
\[
b = b_2 = \frac{h}{\EE((1-p)\rho_0^{\kappa_r}\ln\rho_0)}\, ,
\]
 we see that for this value of $b$, $g_b'(\lambda)/h  = 
1 - \frac{ \EE((1-p)\rho_0^\lambda\ln\rho_0)}{\EE((1-p)\rho_0^{\kappa_r}\ln\rho_0)\EE((1-p)\rho_0^\lambda)}$, and 
$g'_b(\lambda)= 0$ for $\lambda =\kappa_r$. We conclude
\[
 \inf_{b_2 > 0}\Big\{ \sup_{\lambda>0}\Big\{\lambda h- b_2
\ln \EE \big((1-p)\rho_0^\lambda\big)\Big\}\Big\} 
     \leq \kappa_r h,
\]
and so Lemma~\ref{psitildesolution} is proved.
\qed

\begin{lmm}
\label{trapsexist}
Assume that $\EE \ln\rho_0\ne 0$ and let $\kappa$ be defined as after (\ref{not_sinai}).
Let $\gamma > 0$ and fix any $\eps<\gamma/\kappa$. Then, for $\mu$-almost all $\omega$, there is $n_0(\omega, \eps)$ such that for all $n \geq n_0(\omega, \eps)$, 
there is a trap of depth $\Big(\frac{\gamma}{\kappa}-\eps\Big) \ln n$ in the interval $[0, n^\gamma]$.
\end{lmm}

\noindent
{\it Proof.}  
Recall~(\ref{psidef}) and Lemma~\ref{traps}, and keep in mind that the obstacles
are independent from the transition probabilities. We will show that
\begin{equation}
\label{psykh}
\inf\limits_{b_1, b_2} \psi(h, b_1, b_2) =\kappa h\, .
\end{equation}
By~\eqref{eq_sem_traps}, this implies that
\begin{equation}
\label{eq_trapsopt }
\lim_{n\to\infty}\frac{\ln\sup_{b_1, b_2} \PP[\Lambda^x(h, b_1, b_2, n)]}{\ln n}
= -\kappa h\, .
\end{equation}
Take $\eps<\gamma/\kappa$ and chose  $b_1, b_2$ such that for all $n$ large enough
\[
 \PP\Big[\Lambda^x\Big(\frac{\gamma}{\kappa}-\eps, \; b_1, b_2, n\Big)\Big]
\ge n^{-(\gamma-\frac{\kappa \eps}{2})}
\]
(here we use~\eqref{eq_trapsopt }). Divide the interval $[0, n^\gamma]$ into 
$n^\gamma/((b_1+b_2)\ln n)$ intervals of length $(b_1+b_2) \ln n$. Then, 
\begin{align*}
&\PP\Big[\text{there is at least one trap of depth } 
  \Big(\frac{\gamma}{\kappa}-\eps\Big) \ln n \text{ in the interval }[0, n^\gamma]\Big]\\
&~~~~~~~\ge  1-\Big( 1-n^{-(\gamma-\frac{\kappa \eps}{2})}\Big)^{\frac{n^\gamma}{(b_1+b_2)\ln n}}.
\end{align*}
Lemma~\ref{trapsexist} now follows from the Borel-Cantelli lemma.

Now, to prove~\eqref{psykh},  assume that $\EE(\ln\rho_0)< 0$, hence~$\kappa$ is such that 
$\EE(\rho_0^\kappa)=1$ (the case $\EE(\ln\rho_0)> 0$ follows by symmetry).
Using Jensen's inequality, $\ln \EE\big(\frac{1}{ \rho_0^\lambda}\big) \geq \lambda \EE \big(\ln \frac{1}{\rho_0}\big) > 0$, hence we have
\[
 \sup_{\lambda>0}\Big\{\lambda h- b_1\ln \EE \Big(\frac{1}{\rho_0^\lambda}\Big)\Big\}
\leq 
\sup_{\lambda>0}\Big\{\lambda h- \lambda b_1 \EE \Big(\ln \frac{1}{\rho_0}\Big)\Big\},
\]
  and for $b_1 > h \Big(\EE \Big(\ln \frac{1}{\rho_0}\Big)\Big)^{-1}$, 
\[
  \sup_{\lambda>0}\Big\{\lambda h- \lambda b_1 \EE \Big(\ln \frac{1}{\rho_0}\Big)\Big\} =0, 
\]
so we obtain
\[
\inf_{b_1 > 0}\Big\{ \sup_{\lambda>0}\Big\{\lambda h- b_1\ln \EE \Big(\frac{1}{\rho_0^\lambda}\Big)\Big\}\Big\} =0\, .
\]
It remains to show that
\[
\inf_{b_2 > 0}\Big\{ \sup_{\lambda>0}\Big\{\lambda h- b_2\ln \EE (\rho_0^\lambda)\Big\}\Big\} =\kappa h,
\]
but here we can follow verbatim the proof of~\eqref{rightkap} from Lemma \ref{psitildesolution} with $p=0$.
\qed

\medskip

Next, we need to recall some results about hitting and confinement 
(quen\-ched) probabilities for one-dimensional random walks in random environment. 
Obstacles play no role in the rest of this section.
For the proof of these results, see~\cite{CP03} (Sections~3.2 and~3.3)
and~\cite{FGP} (Section~4).

Let $I=[a,c]$ be a finite interval of~$\Z$ with a potential~$V$ defined as in Definition \ref{df_V} and without obstacles.
Let $b$ the first point with minimum potential, i.e.,
\[
b=\min\{x\in [a,c]: \; V(x)=\min_{y\in[a,c]} V(y)\}.
\]
 Let us introduce the following quantities (which depend on~$\om$) 
\begin{align*}
H_-&=\max_{x\in[a, c]} \Bigl(\max_{y \in [a,x]} V(y) -\min_{y\in [x,c]} V(y)\Bigr),\\
H_+&=\max_{x\in[a, c]} \Bigl(\max_{y \in [x,c]} V(y) -\min_{y\in [a,x]} V(y)\Bigr),
\end{align*}
and $H=H_-\wedge H_+$.

First, we need an upper bound on the probability of confinement
in an interval up to a certain time:
\begin{lmm}
\label{l_up_confinement}
There exist $\Upsilon_1,\Upsilon_2>0$ (depending on $\eps_0$), such that for all $u\geq 1$ 
\[
\max_{x\in [a,c]} \Ps^x\Bigl[\frac{T_{\{a,c\}}}{\Upsilon_1(c-a)^4e^H}>u \Bigr] \leq e^{-u},
\]
for $c-a>\Upsilon_2(\eps_0)$.
\end{lmm}

\noindent
{\it Proof.} See Proposition~4.1 of~\cite{FGP}
(in fact, Lemma \ref{l_up_confinement} is a simplified version of that proposition, since
here, due to Condition~E, the potential has bounded increments). \qed

%
%

\medskip

Next, we obtain a lower bound on the confinement probability in the following lemma.
\begin{lmm}
\label{l_low_confinement}
Suppose that $a< b < c$ and that $c$ has maximum potential on $[b,c]$ and $a$ has maximal potential on $[a,b]$.
Then, there exist $\Upsilon_3,\Upsilon_4>0$, such that for all $u \geq 1$ and $x\in (a,c)$
\[
\Ps^x\Bigl[ \Upsilon_3 \ln(2(c-a)) \frac{T_{\{a,c\}}}{e^H}\geq u\Bigr] \geq 
    \frac{1}{2(c-a)}e^{-u},
\]
for $c-a\geq \Upsilon_4$.
\end{lmm}

\noindent
{\it Proof.} See Proposition~4.3 of~\cite{FGP}. \qed

Let us emphasize that the estimates in Lemmas~\ref{l_up_confinement}
and~\ref{l_low_confinement} are valid for all
environments satisfying the uniform ellipticity Condition~E.
An additional remark is due about the usage of Lemma~\ref{l_low_confinement}
(lower bound for the quenched probability of confinement). 
Suppose that there is a trap of depth~$H$ on interval~$[a,c]$,
being~$b$ the point with the lowest value of the potential.
Suppose also that~$a'$ has maximal potential on~$[a,b]$
and~$c'$ has maximal potential on~$[b,c]$. Then, for any~$x\in (a',c')$,
it is straightforward to obtain a lower bound for the probability 
of confinement in the following way: write 
$\Ps^x[T_{\{a,c\}}\geq t]\geq\Ps^x[T_{\{a',c'\}}\geq t]$, and
then use Lemma~\ref{l_low_confinement} for the second term.
This reasoning will usually be left implicit when we use
Lemma~\ref{l_low_confinement} in the rest of this paper.

\subsection{Proofs of Theorems~\ref{t_dim1_ann}--\ref{fixed_omega}.}  
\label{s_proofs_dim1}
{\it Proof of Theorem~\ref{t_dim1_ann}}. 
By Lemma~\ref{traps} and Lemma~\ref{l_low_confinement}, we have for all $b_1, b_2 \in (0, \infty)$ and any $\eps>0$, 
\begin{align*}
\PA[\tau>n]&\ge \PP[A^0(1, b_1, b_2, n)]\\
& \quad  \times \inf_{\s\in A^0( 1, b_1, b_2, n)} \Ps[\xi_t\in (-b_1\ln n, b_2\ln n), \; 
\text{ for all }t\le n ]\\
&\ge  C_1 n^{(b_1+b_2)\ln (1-p)-\psi(1, b_1, b_2)-\eps}
\end{align*}
for all~$n$ large enough.
Thus, recalling that
\[
\widetilde\psi(1)=\inf_{b_1, b_2>0}\{-(b_1+b_2)\ln (1-p)+\psi(1, b_1, b_2)\},
\]
we obtain
\begin{equation}
\label{eq_p22}
\PA[\tau>n]\ge C_2 n^{-\widetilde\psi(1)-\eps}.
\end{equation}

Let us now obtain an upper bound on $\PA[\tau>n]$. Fix $n$, $\beta>0$, $0<\delta<1$. 
We say that the environment $\s$ is good, 
if the maximal (obstacle free) trap depth is less than $(1-\delta)\ln n$ 
in the interval $[-\ln^{1+\beta} n, \ln^{1+\beta} n]$, that is, 
for all $b_1, b_2>0$, $x\in [-\ln^{1+\beta} n, \ln^{1+\beta} n]$ the event
$A^x(1-\delta, b_1, b_2, n)$ does not occur, 
 and also 
\[
\min\{|\Theta\cap[-\ln^{1+\beta} n,0]|, |\Theta\cap[0, \ln^{1+\beta} n]|\}
 \ge \frac{p\ln^{1+\beta} n}{2}.
\]
For any $\eps>0$ we obtain  that for all large enough~$n$ 
\begin{equation}
\label{prob_notgood}
 \PP[\s\text{ is not good}] 
      \le  C_3n^{-\widetilde\psi(1-\delta)+\eps}\ln^\beta n +e^{-C_4 \ln^{1+\beta} n},
\end{equation}
by Lemma~\ref{traps}.

Note that if~$\s$ is good, then 
for every interval $[a,b]\subset [-\ln^{1+\beta} n, \ln^{1+\beta} n]$ 
such that $a,b\in \Theta$, $\Theta\cap (a,b)=\emptyset$ we have 
\[
\max_{x\in[a,b]}V(x)-\min_{x\in[a,b]}V(x)\le (1-\delta)\ln n.
\]
Thus, for such an interval $[a,b]$, on the event $\{ \s\text{ is good}\}$,
 Lemma~\ref{l_up_confinement} (with $u = n^{\delta/2}$) implies that for any $x\in [a, b]$ we have
\begin{equation}
\label{skoko}
\Ps\big[\xi_t\notin [a,b] \text{ for some }t\leq n^{1-\frac{\delta}{2}}
   \mid  \xi_0=x \big]
\ge 1-\exp\Big[-\frac{n^{\delta/2}}{16 \Upsilon_1 \ln^{4+4\beta }n}\Big].
\end{equation}
Let
\[
G=\{\xi_t\notin [-\ln^{1+\beta} n, \ln^{1+\beta} n] \text{ for some }t\le n\}.
\]
Then, by~\eqref{skoko}, on the event $\{ \s\text{ is good}\}$,  we have
(denoting $X$ a random variable with Binomial$\Big(n^{\delta/2}, 1-\exp\Big[-\frac{n^{\delta/2}}{16 \Upsilon_1 \ln^{4+4\beta }n}\Big]\Big)$
distribution) 
\begin{align}
\Ps[\tau>n]&= \Ps[\tau>n, \; G]+\Ps[\tau>n, \; G^c]\nonumber\\ 
&\le (1-r)^{\frac{p}{2}\ln^{1+\beta} n}+(1-r)^{\frac{1}{2}n^{\delta/2}}+ 
\Ps\Big[X\leq \frac{n^{\delta/2}}{2}\Big] \nonumber\\
&\le  e^{ -C_5 \ln^{1+\beta} n}.\label{raz}
\end{align}
To explain the last term in the second line of~~\eqref{raz}, for the event $\{\tau>n\} \cap G^c$, split the time into $n^{\delta/2}$ intervals of length $n^{1-\delta/2}$.
 By~(\ref{skoko}), on each such interval, 
we have a probability of at least $1- \exp\Big[-\frac{n^{\delta/2}}{16 \Upsilon_1 \ln^{4+4\beta }n}\Big]$ to hit an obstacle.
Let~$X'$ count the number of time intervals where that happened. Then, clearly, $X'$ dominates~$X$.

So, by~\eqref{prob_notgood} and~\eqref{raz}, we have
\begin{align*}
 \PA[\tau>n]&\le \int_{\{\s: \; \s\text{ is good}\}}\Ps[\tau>n]\,d\PP
+\int_{\{\s: \; \s\text{ is not good}\}}\Ps[\tau>n]\,d\PP\\
&\le e^{-C_5 \ln^{1+\beta} n}+C_3n^{-\widetilde\psi(1-\delta)+\eps}\ln^\beta n +e^{-C_4 \ln^{1+\beta} n}\\
&\le  C_6n^{-\widetilde\psi(1-\delta)+\eps}\ln^\beta n.
\end{align*}
Together with~\eqref{eq_p22} and Lemma~\ref{psitildesolution},
this implies~\eqref{dim1_ann}.
\qed

\medskip
\noindent
{\it Proof of Theorem~\ref{dim_1_qu}.} 
First we obtain a lower bound on $\Ps[\tau>n]$.
Fix~$a$ and let $b_1, b_2$ be such that 
\[
\widetilde \psi(a)=-(b_1+b_2)\ln (1-p)+\psi(a, b_1, b_2).
\]
Let us show that such $b_1, b_2$ actually exist, that is,
the infimum in~\eqref{psitilddef} is attained. For that, one may reason as follows.
First, since $\psi\geq 0$, for any~$M_0$ there is $M_1>0$ such that 
if $\min\{b_1,b_2\}\geq M_1$ then
\begin{equation}
\label{min_b12}
 -(b_1+b_2)\ln (1-p)+\psi(a, b_1, b_2) \geq M_0.
\end{equation}
Then, it is clear that for any fixed $h>0$ we have
$\lim_{b_2\downarrow 0}\sup_{\lambda>0}(\lambda h-b_2\EE\rho_0^\lambda)=+\infty$,
so (also with the analogous fact for $\frac{1}{\rho_0}$) we obtain that 
for any~$M_0$ there exists $\eps>0$ such that if $\min\{b_1,b_2\}<\eps$ then~\eqref{min_b12}
holds. Thus, one may suppose that in~\eqref{psitilddef} $b_1$ and~$b_2$ vary over a compact set,
and so the infimum is attained.

Let us call such a trap (with $b_1,b_2$ chosen as above) an $a$-optimal trap.
Note that
\[
-\widetilde\psi(a)=\lim_{n\to\infty}\frac{\ln\PP[A^0(a, b_1, b_2, n)]}{\ln n}.
\]
Note also that, since $\widetilde \psi(a)$ is an increasing function 
with $\widetilde \psi(0)=0$ and $1-a$ is a decreasing function, 
there exists unique~$\overline{a}$ such that 
$1-\overline{a}=\widetilde\psi(\overline{a})$.

Fix $\eps' < \frac{\eps}{2}$. Consider an interval $[-n^{\widetilde\psi(\overline{a})+\eps}, n^{\widetilde\psi(\overline{a})+\eps}]$. 
In this interval there will be at least one $\overline{a}$-optimal trap free of 
obstacles with probability at least
\[
1-(1-n^{-\widetilde\psi(\overline{a})-\eps'})^{\frac{n^{\widetilde\psi(\overline{a})+\eps}}{(b_1+b_2)\ln n}}\ge 1-\exp\big(-C_1n^{\eps/2}\big). 
\]
To see this, divide the interval $[-n^{\widetilde\psi(\overline{a})+\eps}, n^{\widetilde\psi(\overline{a})+\eps}]$ into
 disjoint intervals of length $(b_1+b_2)\ln n$ and note that the probability 
that such an interval is an $\overline{a}$-optimal trap free of obstacles, 
by Lemma~\ref{traps}, is at least $n^{-\widetilde\psi(\overline{a})-\eps'}$.
So, a.s.\ for all~$n$ large enough, there  there will be at least one $\overline{a}$-optimal trap free of obstacles in  $[-n^{\widetilde\psi(\overline{a})+\eps}, n^{\widetilde\psi(\overline{a})+\eps}]$. 
If there is at least one $\overline{a}$-optimal trap in  $[-n^{\widetilde\psi(\overline{a})+\eps}, n^{\widetilde\psi(\overline{a})+\eps}]$, 
and $\xi_t$ enters this trap, by Lemma~\ref{l_low_confinement}, 
it will stay there for time~$n^{\overline{a}}$ with probability at least 
$\big(2(b_1+b_2)\ln n\big)^{-(\Upsilon_3+1)}$.
 We obtain
\begin{equation}
\label{q_lower_bound}
\Ps[\tau>n]\ge (r\eps_0)^{n^{\widetilde\psi(\overline{a})+\eps}} 
  \big(2(b_1+b_2)\ln n\big)^{-(\Upsilon_3+1)n^{1-\overline{a}}}
 \ge \exp\Big(-C_2 n^{\widetilde\psi(\overline{a})+\eps}\Big)
\end{equation}
a.s.\ for all~$n$ large enough. 
The factor $(r\eps_0)^{n^{\widetilde\psi(\overline{a})+\eps}} $ appears in~\eqref{q_lower_bound}  because the random walk should first reach the obstacle free $\overline{a}$-optimal trap, 
and on its way in the worst case it could meet obstacles in every site.

To obtain the upper bound, take $\delta>0$ and
consider the time intervals $I_k=[(1+\delta)^k, (1+\delta)^{k+1})$.
If $n\in I_k$, we have
\[
\Ps[\tau>n]\le \Ps[\tau>(1+\delta)^k].
\]

Denote 
\begin{align*}
B_1^{(k)}&=\{\text{the maximal depth of an obstacle free trap in }\\
& \qquad [-(1+\delta)^{k(\widetilde\psi(\overline{a})-\eps)},(1+\delta)^{k(\widetilde\psi(\overline{a})-\eps)}] \text{ is at most } \overline{a} k\ln (1+\delta)\},\\
B_2^{(k)}&=\Big\{|\Theta\cap[-(1+\delta)^{k(\widetilde\psi(\overline{a})-\eps)},0]|
 \ge \frac{p (1+\delta)^{k(\widetilde\psi(\overline{a})-\eps)}}{2}\Big\},\\
B_3^{(k)}&=\Big\{|\Theta\cap[0, (1+\delta)^{k(\widetilde\psi(\overline{a})-\eps)}]|
    \ge \frac{p (1+\delta)^{k(\widetilde\psi(\overline{a})-\eps)}}{2}\Big\}\\
B_4^{(k)}&=\Big\{\text{the maximal length of an obstacle free interval in }\\
& \qquad [-(1+\delta)^{k(\widetilde\psi(\overline{a})-\eps)},(1+\delta)^{k(\widetilde\psi(\overline{a})-\eps)}] \text{ is at most } \frac{2}{\ln\frac{1}{1-p}} k\ln (1+\delta)\Big\}.
\end{align*}

First, we have that for some $C_3$
\[
 \PP[B_i^{(k)}] \geq 1-\exp\big(-C_3 (1+\delta)^{k(\widetilde\psi(\overline{a})-\eps)}\big),
\]
$i=2,3$.
Also,
\begin{align*}
 \PP[B_4^{(k)}] &\geq
 1-2(1+\delta)^{k(\widetilde\psi(\overline{a})-\eps))}
(1-p)^{\frac{2}{\ln\frac{1}{1-p}} k\ln (1+\delta)}\\
&=1-2(1+\delta)^{k(-2+\widetilde\psi(\overline{a})-\eps)}\\
&\geq 1-2(1+\delta)^{-k}.
\end{align*}
Note that, by Lemma~\ref{traps}, in the interval $[-n^{\widetilde\psi(\overline{a})-\eps}, n^{\widetilde\psi(\overline{a})-\eps}]$ 
there will be an obstacle free trap of depth
$\overline{a}\ln n$ with probability  at most 
\[
2n^{\widetilde\psi(\overline{a})-\eps}n^{-\widetilde\psi(\overline{a})+\eps/2}\le C_4 n^{-\eps/2}.
\] 
So, for any $n\in I_k$, in the  interval  $[-(1+\delta)^{k(\widetilde\psi(\overline{a})-\eps)}, (1+\delta)^{k(\widetilde\psi(\overline{a})-\eps})]$ there will be obstacle free traps of depth
$\overline{a}k\ln (1+\delta)$ with probability  at most $C_4 (1+\delta)^{-\eps k/2}$.
Thus, the Borel-Cantelli lemma implies that for almost all~$\s$, a.s.\
for all~$k$ large enough the event $B^{(k)}=B_1^{(k)}\cap B_2^{(k)}\cap B_3^{(k)}\cap B_4^{(k)}$ occurs.  
Analogously to~\eqref{skoko}, for any interval 
$[a,b]\subset [-(1+\delta)^{k(\widetilde\psi(\overline{a})-\eps)}, (1+\delta)^{k(\widetilde\psi(\overline{a})-\eps)}]$ such that 
$(a,b)\cap \Theta=\emptyset$,
we have 
\begin{align}
&\Ps[\xi_t\notin [a,b] \text{ for some }t\leq (1+\delta)^{k(a+\gamma)}
   \mid  \xi_0=x,\; B^{(k)}]\nonumber\\
&~~~~~~~~~~~\ge 1-\exp\Big(-\frac{(1+\delta)^{k\gamma}}{16^2\Upsilon_1\Big(\frac{k\ln(1+\delta)}{\ln(1-p)}\Big)^4}\Big)
\label{skoko1}
\end{align}
for $\gamma>0$ and any $x\in [a,b]$.
Again, let
\[
G^{(k)}=\{\xi_t\notin [-(1+\delta)^{k(\widetilde\psi(\overline{a})-\eps)},(1+\delta)^{k(\widetilde\psi(\overline{a})-\eps)}] \text{ for some }t\le (1+\delta)^k\}.
\]
With the same argument as in \eqref{raz}, we have, for all $\s \in B^{(k)}$
\[
\Ps[\tau>(1+\delta)^k, \; (G^{(k)})^c]
\leq 
(1-r)^{(1+\delta)^{k(1-a-\gamma)}}.
\]
Using \eqref{skoko1}, for $n\in[(1+\delta)^k, (1+\delta)^{k+1})$, we obtain for all $\s \in B^{(k)}$
\begin{align}
\Ps[\tau>n]
&\le \Ps[\tau>(1+\delta)^k ]\nonumber\\
&= \Ps[\tau>(1+\delta)^k, \; G^{(k)}]  
+\Ps[\tau>(1+\delta)^k, \; (G^{(k)})^c]\nonumber\\
&\le (1-r)^{\frac{p  (1+\delta)^{k(\widetilde\psi(\overline{a})-\eps)}}{2}}+(1-r)^{(1+\delta)^{k(1-\overline{a}-\gamma)}}\nonumber\\
&\le e^{-C_3 (1+\delta)^{k(\widetilde\psi(\overline{a})-\eps)}}\nonumber\\
&\le e^{-C_3 n^{\widetilde\psi(\overline{a})-\eps}}.\label{raz1}
\end{align}
Since~$\eps$ is arbitrary and, by Lemma~\ref{psitildesolution},
$\widetilde\psi(\overline{a})=\frac{\kappa_\ell (p) + \kappa_r (p)}
{1+ \kappa_\ell (p) + \kappa_r (p)}$, \eqref{q_lower_bound}
together with~\eqref{raz1} imply~\eqref{qu}.
\qed

\medskip
\noindent
{\it Proof of Theorem~\ref{fixed_obstacles}.}
Denote
\[
u_\beta=\frac{\ln\Big(\frac{1-\beta_0}{\beta_0}\Big)}{\ln\Big(\frac{1-\beta_0}{\beta_0}\Big)
                                          +\ln\Big(\frac{1-\beta_1}{\beta_1}\Big)}
\]
and
\begin{equation}
\label{def_gam}
 \gamma=\frac{|\ln(1-p)|}{|\ln(1-p)|+ F_e},
\end{equation}
where $F_e$ was defined in \eqref{Fedef}. We shall see that $F_e k$ is the maximal possible depth of a trap located at $[a, c]$ with 
$c-a \leq  k $.
Let $B_{n, \alpha}$ be the event that in the interval $[-n^\gamma, n^\gamma]$ there is (at least) one interval of length at least $\alpha\ln n$ which is free of obstacles
and let $I_n(\tb)$ be the biggest such interval. 
For any $\alpha<\gamma|\ln(1-p)|^{-1}$, the event $B_{n, \alpha}$ happens a.s. for $n$ large enough. 
Take an interval $I=[a, c]\subset I_n(\tb) $ and such that $c-a=\lfloor \alpha\ln n\rfloor$. 
For small $\delta$ denote 
\begin{align*}
 U_\delta=\big\{\omega_i^+ & \ge 1-\beta_1 - \delta\text{ for all }i\in [a, a+u_\beta(c-a ],\\
           \omega_i^+ & \le \beta_0+\delta\text{ for all }i\in (a +u_\beta(c-a), c]\big\}.
\end{align*}
So, $U_\delta$ implies that on the interval $[a, c]$ there is a trap of depth $(F_e-\delta')\alpha\ln n$, free of obstacles, where 
$\delta'\to 0$, when $\delta\to 0$. 
Note that
\[
 \mu [U_\delta]\ge K_1(\delta)^{\alpha \ln n}=n^{\alpha \ln K_1(\delta)},
\]
where
\[
 K_1(\delta)=\min\{\mu(\omega_0^+\ge 1-\beta_1-\delta), \mu(\omega_0^+\le \beta_0+\delta)\}.
\]

For $\omega\in U_\delta$, using Lemma~\ref{l_low_confinement}, we obtain
\begin{equation}
 \label{twoterms}
 \Ps[\tau>n]\ge (\eps_0(1-r))^{n^\gamma} \exp\big(-n^{1-\alpha(F_e-2\delta')}\big).
\end{equation}
Now,
\begin{align}
 \Pt[\tau>n] &\geq \mu[U_\delta]  
\int\limits_{\{\omega \in U_\delta\}}\Ps[\tau>n]\mu(d\omega) \nonumber\\
 &\geq n^{\alpha \ln K_1(\delta)}
(\eps_0(1-r))^{n^\gamma}  \exp\big(-n^{1-\alpha(F_e-2\delta')}\big). \label{coisa*}
\end{align}
Note that, by definition~\eqref{def_gam} of $\gamma$, we have $\gamma=1-\gamma F_e |\ln(1-p)|^{-1}$.
Since~$\alpha$ can be taken arbitrarily close to $\gamma |\ln(1-p)|^{-1}$, 
we obtain
\begin{equation}
 \label{coisa**}
\limsup_{n\to\infty}\frac{\ln(-\ln\Pt[\tau>n])}{\ln n}\le \frac{|\ln(1-p)|}{|\ln(1-p)|+F_e}.
\end{equation}

To obtain the other bound, we fix $\alpha>\gamma|\ln(1-p)|^{-1}$. On the event $B_{n, \alpha}^c$, in each of the intervals $[-n^\gamma, 0)$, $[0, n^\gamma]$
there are at least $n^\gamma(\alpha\ln n)^{-1}$ obstacles. Since $F_e k$ is the maximal depth of a trap on an interval of length $k$, the environment $\om$ in the interval  
$[-n^\gamma, n^\gamma]$ satisfies a.s. for large $n$ that the depth of any trap free of obstacles is at most $F_e\alpha \ln n$. Thus, by Lemma~\ref{l_up_confinement}, for any 
$\delta>0$,
the probability that the random walk stays in a trap of depth $F_e\alpha \ln n$
at least for the time $\exp\big((F_e\alpha+\delta) \ln n\big)$ is at most $e^{-C_1 n^{\delta/2} }$. 
We proceed similarly to~\eqref{raz1}. 
Consider the event
\[
 A=\{\xi_t\notin [-n^\gamma, n^\gamma] \text{ for some }t\le n\}.
\]
 We then have
\begin{align*}
\Ps[\tau>n]&=\Ps[\tau>n, A]+\Ps[\tau>n, A^c]\\
&\le (1-r)^{n^\gamma(\alpha\ln n)^{-1}}+((1-r)(1-e^{-C_1 n^{\delta/2}}))^{n^{1-F_e\alpha-\delta}}.
\end{align*}
Since $\delta$ is arbitrary, analogously to the derivation of~\eqref{coisa**} from~\eqref{coisa*} 
one can obtain
\[
\liminf_{n\to\infty}\frac{\ln(-\ln\Pt[\tau>n])}{\ln n}\ge \frac{|\ln(1-p)|}{|\ln(1-p)|+F_e}.
\]
This concludes the proof of Theorem~\ref{fixed_obstacles}.
\qed

\medskip
\noindent 
{\it Proof of Theorem~\ref{fixed_omega}.}
Define 
\begin{align*}
R_n^{\eps, +}(\om)&=\min\{x>0: \; \text{ in $(0,x)$ there is a trap of depth $(1-\eps)\ln n$}\},\\
R_n^{\eps, -}(\om)&=\min\{x>0: \; \text{ in $(-x,0)$ there is a trap of depth $(1-\eps)\ln n$}\},
\end{align*}
and let 
 $R_n^\eps(\om)=\min\{R_n^{\eps,+}(\om), \;R_n^{\eps, -} (\om)\}$. For $\eps=0 $, denote 
$R_n(\om):=R_n^0(\om)$. The statements for $R_n^\eps(\om)$ and $R_n(\om)$ now follow from the definition of these random variables and the invariance principle. 

\medskip
\noindent
{\it Proof of (i).} To obtain a lower bound, we just observe that by Lemma~\ref{l_low_confinement}, when the random walk is in a trap of depth
$\ln n$, it will stay there up to time~$n$ with a probability bounded away from~$0$ (say, by $C_1$).
Further, with $\nu$-probability $(1-p)^{R_n(\om)+1}$ there will be no obstacles in the interval $[0,R_n(\om)]$.
Thus, by the uniform ellipticity,
\[
  \Po[\tau>n]\ge C_1(1-p)^{R_n(\om)+1}\eps_0^{R_n(\om)}\geq e^{-K_2 R_n(\om)}.
\]

To obtain an upper bound, we say that $\tb$ is $k$-good, if
\[
 \big|\Theta\cap [-j, 0]\big|\ge \frac{pj}{2} 
\quad \text{ and }\quad \big|\Theta\cap [0, j]\big|\ge \frac{pj}{2}
\]
for all $j\ge k$. Then,
\begin{equation}
 \label{nu_thetanotgood}
\nu[\tb \text{ is not $R_n^\eps(\omega)$-good}]\leq e^{-C_2 R_n^\eps(\om)}.
\end{equation}
By Lemma~\ref{l_up_confinement}, for all large enough~$n$
\[
 \Ps[\xi_t\in [-R_n^\eps(\om), R_n^\eps(\om)]\text{ for all }t\le n]\le e^{-C_3n^{\eps/2}},
\]
and so on the event $\{\tb\hbox{ is } R_n^\eps(\om)\text{-good}\}$, we have
\begin{align*}
 \Ps[\tau>n]\le e^{-C_3n^{\eps}}+(1-r)^{R_n^\eps(\om)p/2}.
\end{align*}
Together with~\eqref{nu_thetanotgood}, this proves part~(i), since it is elementary to obtain 
that~$R_n^\eps$ is subpolynomial for almost all~$\om$. Indeed, as discussed in Remark~\ref{rem_HuShi},
a comparison to~\cite{HS} suggests that, for~$C_4$ large enough, 
$R_n^\eps(\om)\leq C_4\ln^2n\ln\ln\ln n$ for all but finitely many~$n$, $\mu$-a.s.
Anyhow, one can obtain a weaker result which is still enough for our needs: $R_n^\eps(\om)\leq\ln^4n$
$\mu$-a.s.\ for all but finitely many~$n$, with the following reasoning: any interval of length
$\ln^2 n$ contains a trap of depth at least $\ln n$ with constant probability. So, dividing 
$[-\ln^4n, \ln^4n]$ into subintervals of length~$\ln^2n$, one can see that 
$\mu[\om:R_n^\eps(\om)\geq \ln^4n]\leq e^{-C_5\ln^2n}$ and use the Borel-Cantelli lemma.

\medskip
\noindent
{\it Proof of (ii).} Take  $a=\frac{\kappa}{\kappa+1}$ and $0<\eps<a/\kappa$. In this case, by Lemma~\ref{trapsexist} 
there is a trap of depth $\left(\frac{a}{\kappa}-\eps\right)\ln n$ on the interval $[0, n^a)$, a.s.\ for all~$n$ large enough.
Using Lemma~\ref{l_low_confinement}, when the random walk enters this trap, 
it stays there up to time~$n$ with probability at least
$\frac{1}{2n^a}\exp\big(-\Upsilon_3 \ln(2n^a)n^{1-\frac{a}{\kappa}+\eps}\big)$. 
Further, the probability that the interval $[0, n^a)$ is free of obstacles 
is $(1-p)^{n^ a}$, and we obtain
\begin{align*}
 \Po[\tau>n] &\ge (1-p)^{n^ a}\eps_0^{n^ a} 
 \frac{1}{2n^a}\exp\big(-\Upsilon_3 \ln(2n^a)n^{1-\frac{a}{\kappa}+\eps}\big)\\
          &\ge e^{-C_6n^{a+2\eps}}\\
          &=e^{-C_6n^{ \frac{\kappa}{\kappa+1}+2\eps}},
\end{align*}
as~$a$ was chosen in such a way that $a=1-\frac{a}{\kappa}=\frac{\kappa}{\kappa+1}$. 
Since $\eps > 0$ was arbitrary, this yields
\[
\liminf_{n\to\infty}\frac{\ln(-\ln \Po[\tau>n])}{\ln n} \geq\frac{\kappa}{\kappa+1}\, .
\]
To prove the corresponding upper bound, we proceed analogously to the proof of~\eqref{dim1_ann}. We say that the obstacle environment $\tb$ is good (for fixed $n$), if 
\[
\min\{|\Theta\cap[-n^a,0]|, |\Theta\cap[0, n^a]|\}
 \ge \frac{pn^a}{2}.
\]
Note that
\[
\nu[\tb\text{ is good}]\ge 1-e^{-C_7 n^a}.
\]
Let 
\[
 G=\{\xi_t\notin [-n^a, n^a] \text{ for some } t\le n\}.
\]
Observe that, by Theorem 1.2 from~\cite{FGP}, we have $\Po[G]\leq e^{-C_8 n^{1-\frac{a}{\kappa}}}$.
Again, since for $a=\frac{\kappa}{1+\kappa}$ we have $a=1-\frac{a}{\kappa}$,
 on the event $\{\tb \text{ is good}\}$ we obtain
\begin{align*}
 \Po[\tau>n]& \leq \Po[G]+\Po[\tau>n, G^c]\\
&\leq e^{-C_8 n^{1-\frac{a}{\kappa}}}+(1-r)^{\frac{pn^a}{2}}\\
&\leq e^{-C_9n^{ \frac{\kappa}{1+\kappa}}}.
\end{align*}
This concludes the proof of Theorem~\ref{fixed_omega}.
\qed

\section*{Acknowledgements}
S.P.\ and M.V.\ are grateful to Fapesp (thematic grant 04/07276--2), 
CNPq (grants 300328/2005--2, 304561/2006--1, 471925/2006--3) and 
CAPES/DAAD (Probral) for partial support. 


We thank Francis Comets for instructive discussions. Also, we thank the anonymous referees 
for their careful lecture of the first version, which lead to many improvements.

\end{document}